\documentclass[aps, prd, floats, floatfix, onecolumn, superscriptaddress, nofootinbib, letterpaper, notitlepage, preprintnumbers]{revtex4-1}

\usepackage{graphicx} 
\graphicspath{{images/}}
\usepackage{amsmath}
\usepackage{amsfonts}
\usepackage{amssymb}
\usepackage{amstext}
\usepackage{tensor}
\usepackage{fullpage}
\usepackage[english]{babel}
\usepackage{blindtext}
\usepackage{helvet}
\usepackage{microtype}
\usepackage[pdftex,
  pdftitle={Spherical Bessel Functions},
  pdfauthor={Teboho A. Moloi},
  bookmarks,bookmarksopen=false,
  pdfstartview={FitH},
  linktocpage=true
]{hyperref}

\newcommand{\li}{\mathrm{Li}}

\begin{document}

\title{Spherical Bessel Functions}

\author{Teboho A. Moloi}
\email{teboho.abram.moloi@gmail.com}
\affiliation{Department of Physics, Nelson Mandela University, Port Elizabeth, 6031, South Africa}

\date{\today}

\preprint{MIT-CTP/4872}

\begin{abstract}
We examine indefinite integral involving of arbitrary power $x$, multiplied by three spherical Bessel functions of the first kind $j_{h},j_{k}$, and $j_{l}$ with integer order $h,k,l \geq 0$ and an exponential. Then we add some conditions for smooth calculation in considering the general and elementary exact evaluation. Thus, by measuring their equivalents, we can verify their accuracy.
\end{abstract}

\maketitle

\section{Introduction}
In the past, numerous studies have included integrating the zero to four spherical function
of Bessel as form of applications from different fields such as astrophysics, particle physics, nuclear physics etc. \cite{Mehrem:2010qk,Mehrem1991}, had a way of measuring these integrals using a
triangle in an effort to achieve precise numerical integrations of these combinations.
And a very interesting method was also developed involving an exponential, two Bessel and a polynomial equation introduced by \cite{Gebremariam:2010zz}. All of these methods seem restrictive, and \cite{Fabrikant:2013fv}, then presented an alternative solution with no restrictions. Where the relation between spherical function and function of Bessel were taken up and the trigonometric identities was applied to the integral. They focused on the estimation of the infinite integrals that involve polynomial multiplied by three spherical functions of the Bessel and an exponent.
In this text, we have followed this approach and ways of applying such integral conditions to promote numerical analysis. 
\par
Specifically, we focus on indefinite integral of the form
\begin{align}\label{eq:formbess}
I^{nm}_{hkl}(x;\alpha,\beta,\mu) = \int dx x^{n}e^{-mx}j_{h}(\alpha x)j_{k}(\beta x)j_{l}(\mu x).
\end{align}
Here $j_{h},j_{k}$, and $j_{l}$ are spherical Bessel functions of the first kind, this integral converge when $h,k,l\geq0$, and $\alpha,\beta,\mu$ are real numbers.
\par
This paper has the following structure; we provide the general literature on Bessel function in \ref{sec:lit}. In section.~\ref{sec:prelim} we introduce some of the important equations required in our results. We provide some results in section.~\ref{sec:result}. We finally, conclude in section.\ref{sec:concl}.
\section{Spherical Bessel Functions literature}\label{sec:lit}

Spherical Bessel functions are well-known to account for the problems with circular symmetry. In spherical coordinate, if one solve Helmholtz's and Laplacian's equation the solution yield the following differential equation
\begin{align}\label{eq:differen}
\frac{d^2 y}{dx^{2}} + \frac{2}{x}\frac{dy}{dx} + \bigg[1 +\frac{l(l+1)}{x^{2}} \bigg]y = 0.
\end{align}
Spherical Bessel functions with indices that are not integers are usually less important to implement, here we presume that index $l$ to be integral. The solution to expression above \ref{eq:differen} results in the spherical Bessel- and Neumann-function, $j_{l}(kr)$ and $n_{l}(kr)$ they are both respectively defined as follows;
\begin{align}
j_{l}(x) = \sqrt{\frac{\pi}{2x}} J_{l+1/2}(x), \qquad n_{l}(x) = \sqrt{\frac{\pi}{2x}} N_{l+1/2}(x) = (-1)^{l+1}\sqrt{\frac{\pi}{2x}}J_{-l-1/2}(x).
\end{align}
Here we make note that $j_{l}(x)$ are Bessel functions of the first kind and $n_{l}(x)$ are the second kind, we then can relate the first and second kind by the following expression
\begin{align}
h^{(1)}_{l}(x) = \sqrt{\frac{\pi}{2x}} H^{(1)}_{l+1/2}(x) = j_{l}(x) + in_{l}(x),\qquad h^{(2)}_{l}(x) = \sqrt{\frac{\pi}{2x}} H^{(2)}_{l+1/2}(x) = j_{l}(x) - in_{l}(x),
\end{align}
with $h^{(1)}_{l}(x)$  and $h^{(2)}_{l}(x)$ are Spherical Hankel functions and their counter-parts are Hankel functions $H^{(1)}_{l}$ and $H^{(2)}_{l}$, now as we can see from the information provided above $j_{l}$ and $n_{l}$ are spherical Bessel functions and their counter-parts are Bessel functions. The Bessel functions and spherical Bessel functions are related---this can be shown by the function $\sqrt{x}j_{l}(x)$ and $\sqrt{x}n_{l}(x)$ both satisfy the Bessel functions. From the series solution, with the conventional normalization (for more details see \cite{arfken2013mathematical}) one can show that
\begin{align}
y_{l-1} + y_{l+1} = \frac{2l+1}{x}, \qquad ly_{l-1}-(l+1)y_{l+1} = (2l+1) \frac{dy_{1}}{dx},\\
\frac{d}{dx}\bigg[x^{l+1}y_{l}(x)\bigg] = x^{l+1}y_{l-1}(x), \qquad \frac{d}{dx}\bigg[x^{-1}y_{l}(x)\bigg] = x^{-1}y_{l+1}(x),
\end{align}
where $y_{l}$ can be any of the following functions $j_{1}, \, n_{l}, \, h^{(2)}_{l},$ and $h^{(2)}_{l}$. 
This recurrence relations in turn leads back to the differential equations. The spherical Bessel function can be computed by indiction on $l$ which leads to Rayleigh's formulas;
\begin{align}\label{eq:Ray}
j_{l}(x) &= (-1)(-x)^{l}\bigg(\frac{1}{x}\frac{d}{dx}\bigg)\frac{\sin(x)}{x},\qquad n_{l}(x) = (-1)(-x)^{l}\bigg(\frac{1}{x}\frac{d}{dx}\bigg)\frac{\cos(x)}{x},\\
h^{(1)}_{l}(x) &= (-1)(-x)^{l}\bigg(\frac{1}{x}\frac{d}{dx}\bigg)\frac{i}{x}e^{ix}, \qquad h^{(2)}_{l}(x) = (-1)(-x)^{l}\bigg(\frac{1}{x}\frac{d}{dx}\bigg)\frac{i}{x}e^{-ix}.
\end{align}
In this paper, we focus on spherical Bessel functions of the first, second and third kind for integer $l \geq 0$. From the formulas given above, for $l = 0,1,2$ the solutions are denoted by
\begin{align}
j_{0} &= \frac{\sin(x)}{x}, \qquad j_{1} = \frac{\sin(x)}{x^{2}} - \frac{\cos(x)}{x}, \qquad j_{2} = \bigg(\frac{3}{x^{3}}-\frac{1}{x} \bigg)\sin(x) - \frac{3}{x^{2}}\cos(x),\\
n_{0} &= \frac{\cos(x)}{x}, \qquad j_{1} = -\frac{\cos(x)}{x^{2}} - \frac{\sin(x)}{x}, \qquad j_{2} = \bigg(\frac{3}{x^{3}}-\frac{1}{x} \bigg)\cos(x) - \frac{3}{x^{2}}\sin(x),\\
h^{(1)}_{0} & = -\frac{i}{x}e^{ix}, \qquad h^{(1)}_{1} = e^{ix}\bigg(\frac{1}{x} - \frac{i}{x^{2}} \bigg),\qquad h^{(1)}_{2} = e^{ix}\bigg(\frac{i}{x} - \frac{3}{x^{2}} -\frac{3i}{x^{3}}\bigg).
\end{align}
Furthermore, from the Rayleigh expressions provided in \eqref{eq:Ray} we can easily extract limiting behaviors; for examples for $x \ll l$, the solution have the following behavior \cite{cahill_2019}
\begin{align}
j_{l} \approx \frac{2^{l}l!}{(2l+1)!}x^{l} = \frac{x^{l}}{(2l+1)!!}, \qquad n_{l}\sim -\frac{(2l)!}{2^{l}l!}x^{-(l+1)} = -\frac{(2l-1)!!}{x^{l+1}},   
\end{align}
and for $x \gg l$ we then have
\begin{align}
j_{l} \approx \frac{1}{x}\sin\bigg(x-\frac{l\pi}{2}\bigg), n_{l} \approx -\frac{1}{x}\cos\bigg(x-\frac{l\pi}{2}\bigg), \qquad h^{(1)}_{l} \sim (-1)^{1+l}\frac{e^{ix}}{x}.
\end{align}
Now lets introduce the closure relation of the spherical Bessel functions, starting with the non-trivial formula we get 
\begin{align}\label{eq:expansion}
e^{ikz} = \sum^{\infty}_{l=0}(2l+1)i^{l}j_{l}(kr)P_{l}(\cos \theta). 
\end{align} 
where $P_{l}$ is the Legendre polynomial of order $l$, if the wave vector is pointing at the direction than the positive z-axis, then the above expression \eqref{eq:expansion} can be generalized; we make a note that $Y^{0}_{l}(\theta,\phi) =\sqrt{(2l+1)/4\pi} P_{l}(\cos\theta)$, we find
\begin{align}\label{eq:planewave}
e^{i\vec{k}\cdot \vec{x}} = 4\pi\sum^{\infty}_{l=0}i^l j_{1}(kr) \sum ^{l}_{m=-l}Y^{m\ast}(\theta_{\vec{k}}\phi_{\vec{k}})Y^{m}(\theta_{\vec{k}}\phi_{\vec{k}}).
\end{align}
Lets now normalize the delta function, the usefulness of this will be seen later as consequence of the identity of \eqref{eq:planewave} is the inner-product of the two spherical Bessel functions. Solving this we begin with the following
\begin{align}
\int d\vec{x}\,e^{i\vec{k}\cdot\vec{x}}e^{-i\vec{k}'\cdot\vec{x}}  =(2\pi)^{3}\delta(\vec{k} - \vec{k}'),
\end{align}
following from \eqref{eq:planewave} the right hand side yields
\begin{align}
\int d\vec{x}\,e^{i\vec{k}\cdot\vec{x}}e^{-i\vec{k}'\cdot\vec{x}} = \sum_{l,m}(4\pi)^{2}\int dr\,r^{2} j_{l}(kr)j_{1}(k'r)Y^{m\ast}(\theta_{\vec{k}})Y^{m}(\theta_{\vec{k}}).
\end{align}
While on the right side, we get the following
\begin{align}
(2\pi)^{3}\delta(\vec{k} - \vec{k}') = (2\pi)^{3}\frac{1}{k^{2}\sin\theta}\delta(\vec{k} - \vec{k}')\delta(\vec{\theta} - \vec{\theta}')\delta(\vec{\phi} - \vec{\phi}'),
\end{align}
combining the above information and making note of the fact that $\sum_{l,m} Y_{l}^{m\ast}(\Omega_{\vec{k}})Y^{m}_{l}(\vec{k}) = \delta( \Omega_{\vec{k}} - \Omega_{\vec{k}}')$, we arrive at
\begin{align}
\int^{\infty}_{0}dr\,r^{2} j_{l}(kr)j_{1}(k'r) = \frac{\pi}{2k^{2}}\delta(\vec{k}-\vec{k}'),
\end{align} 
with some mathematics we also consider orthogonality relation which reads
\begin{equation}
\int^{\infty}_{-\infty}j_{k}(x)j_{l}(x)dx = \frac{\pi}{2l+1} \delta_{kl}, \quad \text{for} \quad {k,l \in  \mathbb{N}}, 
\end{equation}
here $\delta_{kl}$ is Kronecker delta. It is useful to mention that infinite integrals over one, two and three Bessel functions over the years have gained interest and they are well-known \cite{Mehrem2011ThePW,Auluck2010SomeII,Mehrem1991,Maxinon1991,Rosenheinrich2016,Qi2016OnTF,Yacsar2016UnifiedBM,Bloomfield2017IndefiniteIO,Adkins2013ThreedimensionalFT,Valery2013}. The integrals are generally expressed in accordance with the prescribed functions of the Bessel functions and some coefficients of these functions are established from a finite series, the terms of which are obtained from recurrence relation which involve polynomials.
\par
Below we introduce the two most important recursion relation which are always fulfilled by spherical Bessel functions--- which normally act as a connector between contiguous $l$, as well as the derivatives $\partial_{x}y_{l}$ and various $y_{l}$: they require a bit of trail and error, which are given by;
\begin{align}
y_{l} &= \pm \partial_{x}y_{l-1}(x) + \frac{l\pm1}{x}y_{l\pm1}(x)\\
y_{l} &=\pm y_{l\pm2}(x) +\frac{2l\pm1}{x}y_{l\pm1}{x}.  
\end{align}
\section{Preliminaries} \label{sec:prelim}
In addition to work presented in ref.~\cite{Bloomfield2017IndefiniteIO} we provide the exponential-integral which take a form
\begin{align}\label{eq:expon}
Z_n(x) = \int dx x^{n}e^{x},
\end{align}
where $n$ denotes an integer. Evaluating this yields a recursion relation, which can be obtained by performing integration by parts in Eq.~\eqref{eq:expon}, which will output the following expression
\begin{align}\label{eq:recursionr}
Z_n(x) = x^{n}e^{x} - nZ_{n-1}(x).
\end{align}
It follows that the recursion will allow us to step down by $n$ for integer values above zero and the is some kind of cut--off at integers strictly equal to zero, and there will be a step up of $n$ for integer values below zero and the is some kind of cut--off at integers strictly equal to negative one \cite{Bloomfield2017IndefiniteIO}, leading to definite integrals which where introduced by Sch$\ddot{o}$milch and Arndt called the sine-integral, cosine-integral shown in ref.~\cite{Bloomfield2017IndefiniteIO} and the exponential-integral \cite{2019arXiv190712373M} with the following form
\begin{align}\label{eq:ei}
Ei(x) = \int^{-x}_{\infty} \frac{e^{-u}}{u}du.
\end{align}
Where $Ei(x)$ is a special function, this exponential-integral was introduced evaluated for all real values  by Sch$\ddot{o}$milch which is related to Logarithm-integral in this manner
\begin{align}\label{eq:relationle}
\li(x) = \int\frac{du}{\log u}, \quad \text{relation},\quad \li\,e^{x} = Ei(x),
\end{align}
We provide a schematic for of the special function $Ei(x)$,
\begin{figure}[ht!]
  \centering
  \includegraphics[scale=0.3]{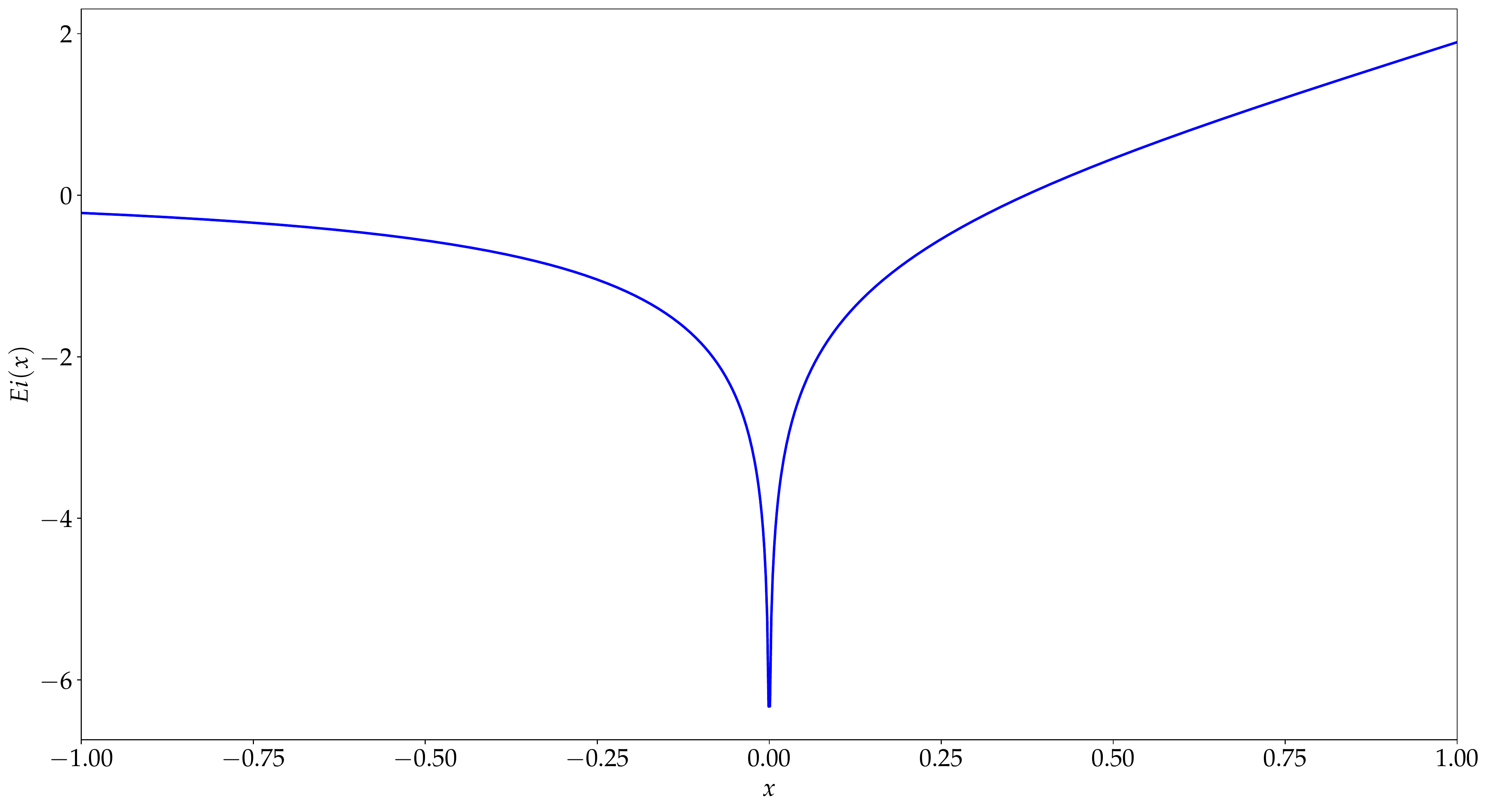}
  \caption{Exponential-integral evaluated at specific range.}\label{fig:exp}
\end{figure}
We see from Fig.~\ref{fig:exp} that the non-conical function dubbed the exponential function, in scale $x<0$ it appears to be negative, it is positive on $x>0$, and it has a unique zero. We also notice that it is concave on scales from negative infinity to zero and alsonononon scales $[0,1]$ and decreases on scales from negative infinity to zero. While it appears to be convex from scales from one to positive infinity and it is increasing on scales from zero to positive infinity.

\section{Integrals of Three Spherical Bessels Of Different Order, Exponent, And Polynomial}\label{sec:result}
We focus on the integral of the form
\begin{align}\label{eq:mainint}
I^{nm}_{hkl}(x;\alpha,\beta, \mu) = \int dx x^{n}e^{-mx}j_{h}(\alpha x)j_{k}(\beta x)j_{l}(\mu x), 
\end{align}
with $n,m \in \mathbb{Z}$, $h,k,l \in \mathbb{N}$, and $\alpha,\beta,\mu \in \mathbb{R}$. These form of Bessel functions are previously, robustly investigated in \cite{Gebremariam:2010zz,Mehrem:2010qk,Mehrem:2011ah,2019arXiv190712373M}. It is easy to see that if one consider the scenario where $\alpha=0$, or  $\beta = 0$, or $\mu = 0$ the integral Eq.~\eqref{eq:mainint} become finite, so we look at scenarios where $\alpha\neq0$, $\beta \neq 0$, and $\mu \neq 0$. The general solution of the integral become \cite{Fabrikant:2013fv}  
\begin{align}\label{eq:genlint}
I^{nm}_{hkl}(x;\alpha,\beta, \mu) &= \frac{\pi}{\alpha\beta\mu \Gamma(3-n)\sin(\pi n)}\bigg\{-[m^{2}+(\alpha+\beta+\mu)^{2}]^{\frac{(2-n)}{2}}\sin\bigg[(n-2)\arctan\bigg(\frac{(\alpha+\beta+\mu)}{m}\bigg)\bigg] \nonumber \\
&+ [m^{2}+(\mu+\beta-\alpha)^{2}]^{\frac{(2-n)}{2}}\sin\bigg[(n-2)\arctan\bigg(\frac{(\mu+\beta+\alpha)}{m}\bigg)\bigg]  \nonumber \\
&+ [m^{2}+(\mu+\alpha-\beta)^{2}]^{\frac{(2-n)}{2}}\sin\bigg[(n-2)\arctan\bigg(\frac{(\mu+\alpha-\beta)}{m}\bigg)\bigg] \nonumber \\ &+
[m^{2}+(\alpha+\beta-\mu)^{2}]^{\frac{(2-n)}{2}}\sin\bigg[(n-2)\arctan\bigg(\frac{(\alpha+\beta-\mu)}{m}\bigg)\bigg]\bigg\}+c
\end{align}
Here $\Gamma$ is a Gamma function. We then extend work by \cite{Fabrikant:2013fv}  by looking at different simple cases, such as when $h,k,l$ are all set to zero, or if $m=i$ or if $m=0$ and so forth.
Now if $h=k=l=0$ we get the following
\begin{align}
I^{nm}_{000}(x;\alpha,\beta, \mu) = \int dx x^{n} e^{-mx}j_{0}(\alpha x)j_{0}(\beta x)j_{0}(\mu x).
\end{align}
Taking into account the trigonometric identity below
\begin{align}
\sin(\alpha x)\sin(\beta x)\sin(\mu x) =-\sin[x(\alpha+\beta+\mu)]+\sin[x(\alpha-\beta+\mu)]+\sin[x(\alpha+\beta-\mu)].
\end{align}
Putting all information together we arrive at
\begin{align}
I^{nm}_{000}(x;\alpha,\beta, \mu) =  -\frac{i}{2\alpha\beta\mu}[\mathcal{A}(\alpha,\beta, \mu)+\mathcal{B}(\alpha,\beta, \mu) +\mathcal{C}(\alpha,\beta, \mu) +c,
\end{align}
where $\mathcal{A},\mathcal{B}$, and $\mathcal{C}$ are defined as follows
\begin{align}
\mathcal{A}&=x^{n+1}[ix(\alpha-\beta+\mu-im)]^{-n-1}\Gamma[n+1,i(\alpha-\beta+\mu-im)x]- \nonumber\\
&x^{n+1}[x(m-i(\alpha-\beta+\mu))]^{-n-1}\Gamma[n+1,(m-i(\alpha-\beta+\mu))x],\label{eq:ex1}\\
\mathcal{B}&= x^{n+1}[ix(\alpha+\beta-\mu-im)]^{-n-1}\Gamma[n+1,i(\alpha+\beta-\mu-im)x]-\nonumber\\
&x^{n+1}[x(m-i(\alpha+\beta-\mu))]^{-n-1}\Gamma[n+1,(m-i(\alpha+\beta-\mu))x],\label{eq:ex2}\\
\mathcal{C}&=-\bigg\{x^{n+1}[ix(\alpha+\beta+\mu-im)]^{-n-1}\Gamma[n+1,i(\alpha+\beta+\mu-im)x]-\nonumber\\
&x^{n+1}[x(m-i(\alpha+\beta+\mu))]^{-n-1}\Gamma[n+1,(m-i(\alpha+\beta+\mu))x]\bigg\}\label{eq:ex3}.
\end{align}
We can express the above functions in the different form, such that
\begin{align}
\mathcal{A}=x^{n+1}\{[E_{\_ n}i(\alpha-\beta+\mu)x]-E_{\_ n}[(m-i(\alpha-\beta+\mu))x]
\end{align}
Now provided Eqs.~\eqref{eq:ex1}---\eqref{eq:ex3} we notice that terms with $im$ if we consider a case we $m=i$ reduce to one, while for cases where $m=0$ all the terms with $m$ vanishes.
\section{Conclusions}\label{sec:concl}
In this paper, we revisited the indefinite integral of the power of $x$ which we multiplied with three spherical Bessel function and an exponential. In this study, we adopted the general and elementary method which introduced by \cite{Fabrikant:2013fv}, we extended this work by providing some limitation on spherical Bessel functions of the first kind with different order. This expressions allow for accurate computations, where the general and elementary method might not be computed smoothly.

\acknowledgments

This work is supported in part by the Nelson Mandela University. Department of Physics under grant council, and in part by NMU Research Development.

\appendix




\bibliographystyle{utphys}
\bibliography{bessel}

\end{document}